\documentclass[fleqn]{article}
\usepackage{amssymb}
\usepackage{amsfonts}
\usepackage{geometry}

\begin{document}

\title{Volume estimates and their applications in the problem of optimal
recovery}
\author{Alexander Kushpel \\
Department of Mathematics, \c{C}ankaya University\\
Ankara, Turkey\\
E-mail: kushpel@cankaya.edu.tr}
\maketitle

\begin{abstract}
We study volumes of sections of convex origin-symmetric bodies in $\mathbb{R}%
^{n}$ induced by orthonormal systems on probability spaces. The approach is
based on volume estimates of \ John-L\"{o}wner ellipsoids and expectations
of norms induced by the respective systems. The estimates obtained allows us
to establish lower bounds for the radii of sections which gives lower bounds
for Gelfand widths (or linear cowidths). As an application we offer a new
method of evaluation of Gelfand and Kolmogorov widths of multiplier
operators. In particular, we establish sharp orders of widths of standard
Sobolev classes $W_{p}^{\gamma }$ in $L_{q}$ in the difficult case, i.e. $%
1<q\leq p\leq \infty $.
\end{abstract}

{\bf Keywords}
Volume, Convex body, Recovery

{\bf Subject}
MSC 41A46, 42C10

\section{Introduction}

In this article we study the problem of optimal recovery of function classes
on general probability spaces $\left( \Omega ,\mathcal{A},\nu \right) $
using linear information. We start with the respective extremal problems.
Consider cowidth of a set $A$ in a Banach space $X$. More precisely, let $%
\left( X,\left\Vert \cdot \right\Vert \right) $ be a Banach space with the
unit ball $B=B_{X}=\left\{ x\left\vert x\in X, \left\Vert
x\right\Vert \leq 1\right. \right\} $ and $A\subset X$. Let $Y$ be a coding
set, i.e. the set which contains information on the elements of $A$ and $%
\Phi :A\rightarrow Y$ be a family of coding operators $\phi :A\rightarrow Y$%
. Let $D\subset X$,%
\[
\mathrm{diam}\left( D\right) =\sup \left\{ \left\Vert x-y\right\Vert
\left\vert x,y\in D\right. \right\}
\]%
and
\[
\phi ^{-1}\left( x\right) =\left\{ y\left\vert y\in X, \phi \left(
y\right) =\phi \left( x\right) \right. \right\}
\]%
be the preimage of $x\in X$. The respective cowidth is defined as
\[
\mathrm{co}^{\Phi }\left( A,X\right) =\inf_{\phi \in \Phi }\sup_{x\in A}%
\mathrm{diam}\left\{ \phi ^{-1}\left( \phi \left( x\right) \right) \right\}
.
\]%
In particular, if $Y=\mathbb{R}^{n}$ and $\Phi =\mathcal{L}\left( \mathrm{lin%
}\left( A\right) ,\mathbb{R}^{n}\right) $ is the set of linear operators
from $\mathrm{lin}\left( A\right) $ into $\mathbb{R}^{n}$ then $\mathrm{co}%
^{\Phi }\left( A,X\right) $ is the linear cowidth $\lambda ^{n}\left(
A,X\right) $. Let $u\in \mathcal{L}\left( X,Z\right) $ be a linear operator,
$u:X\rightarrow Z$, then
\[
\lambda ^{n}\left( u\right) =\mathrm{co}^{\Phi }\left( u\right) =\mathrm{co}%
^{\Phi }\left( uB_{X},Z\right) .
\]%
It is easy to see that $\lambda ^{n}\left( A,X\right) =2d^{n}\left(
A,X\right) $, where%
\[
d^{n}\left( A,X\right) =\inf_{L^{n}}\sup_{x\in L^{n}}\left\Vert x\right\Vert
, \textrm{codim}L^{n}=n
\]%
is the Gelfand width. Let $u\in \mathcal{L}\left( X,Z\right) $, and $u^{\ast
}$ be its adjoint operator. It is well-known that if $u$ is compact or $Z$
is reflexive (see. e.g. \cite{mak}) then%
\begin{equation}
d^{n}\left( u^{\ast }\right) =d_{n}\left( u\right) ,  \label{dual}
\end{equation}%
where%
\[
d_{n}\left( u\right) =\inf_{L_{n}\subset X}\sup_{x\in B_{X}}\inf_{y\in
L_{n}}\left\Vert ux-y\right\Vert _{Z}, \dim L_{n}=n
\]%
is the Kolmogorov width of $uB_{X}$ in $Z$. Similarly, let $A$ be an
origin-symmetric set in $Z$. Kolmogorov width $d_{n}\left( A,Z\right) $ is
defined as%
\[
d_{n}\left( A,Z\right) =\inf_{L_{n}\subset X}\sup_{x\in A}\inf_{y\in
L_{n}}\left\Vert x-y\right\Vert _{Z}, \dim L_{n}=n.
\]

Consider a general probability space $\left( \Omega ,\mathcal{A},\nu \right)
$. Let $L_{p}\left( \Omega ,\mathcal{A},\nu \right) =L_{p}$ be the set of
functions $f:\Omega \rightarrow \mathbb{R}$ of finite norm $\left\Vert
f\right\Vert _{p}$ given by%
\[
\left\Vert f\right\Vert _{p}=\left\{
\begin{array}{cc}
\left( \int_{\Omega }\left\vert f\right\vert ^{p}d\nu \right) ^{\frac{1}{p}},
& 1\leq p<\infty , \\
\mathrm{ess\sup }\left\vert f\right\vert , & p=\infty%
\end{array}%
\right.
\]%
and $U_{p}=\left\{ f\left\vert \left\Vert f\right\Vert _{p}\leq 1\right.
\right\} $ be the unit ball in $L_{p}$. The inner product of functions $f$
and $g$ in $L_{2}$ is defined as usual, $\left[ f,g\right] =\int_{\Omega
}fgd\nu $. Let $\phi _{k}:\Omega \mapsto \mathbb{R}$, $\phi _{k}\in
L_{\infty }$, $k\in \mathbb{N}$, be a fixed orthonormal system on $\left(
\Omega ,\mathcal{A},\nu \right) $. Then the formal Fourier expansion for $%
\varphi \in L_{p}$, $1\leq p\leq \infty $ can be written as
\[
\varphi \sim \sum_{k=1}^{\infty }c_{k}\left( \varphi \right) \phi _{k}
, c_{k}\left( \varphi \right) =\left[ \varphi ,\phi _{k}\right] .
\]%
The Fourier sum $S_{m}\left( \varphi \right) $ of order $m$ is%
\[
S_{m}\left( \varphi \right) =\sum_{k=1}^{m}c_{k}\left( \varphi \right) \phi
_{k}.
\]%
Let $\Lambda =\left\{ \lambda _{k}, k\in \mathbb{N}\right\} $ be a
sequence of real numbers. If for any $\varphi \in L_{p}$ there is a function
$f=\Lambda \varphi \in L_{q}$ such that
\[
f=\Lambda \varphi \sim \sum_{k=1}^{\infty }\lambda _{k}c_{k}\left( \varphi
\right) \phi _{k}
\]%
then we say that the multiplier operator is of $\left( p,q\right) $-type. We
say that $f\in \Lambda U_{p}$ if%
\[
\Lambda \varphi =f\sim \sum_{k=1}^{\infty }\lambda _{k}c_{k}\left( \varphi
\right) \phi _{k},
\]%
where $\varphi \in U_{p}$. Clearly, if $\left\{ \left\vert \lambda
_{k}\right\vert , k\in \mathbb{N}\right\} $ is a non increasing
sequence then
\begin{equation}
\sup \left\{ \left\Vert f-S_{m}\left( f\right) \right\Vert _{2}\left\vert
f\in \Lambda U_{2}\right. \right\} =\left\vert \lambda _{m+1}\right\vert .
\label{FOURIER}
\end{equation}

We study lower bounds of linear cowidths $\lambda ^{n}\left( \Lambda
U_{p},L_{q}\right) $ of function classes $\Lambda U_{p}$ generated by
multiplier orators $\Lambda $ in $L_{q}$. Various examples of such function
classes are presented in Section \ref{EXAMPLES}. The problem of optimal
reconstruction of functions on general domains $\Omega $ in $\mathbb{R}^{d}$
has recently attracted a lot of attention \cite{novak-1}. General bounded
Lipschitz domains were studied in \cite{novak-2}, \cite{triebel-1}, \cite%
{vibiral-1} \cite{Dahlke}, \cite{Mieth-1}.

Let $\mathbf{e}_{1},\cdots ,\mathbf{e}_{n}$ be the canonical basis in $%
\mathbb{R}^{n}$, $\alpha =(\alpha _{1},\cdots ,\alpha _{n})\in \mathbb{R}%
^{n} $, $\beta =(\beta _{1},\cdots ,\beta _{n})\in \mathbb{R}^{n}$ and $%
\left( \alpha ,\beta \right) =\sum_{k=1}^{n}\alpha _{k}\beta _{k}$. Also,
let $\Vert \alpha \Vert _{(2)}=\left( \alpha ,\alpha \right) ^{\frac{1}{2}}$
be the Euclidean norm on $\mathbb{R}^{n}$, $l_{2}^{n}=\left( \mathbb{R}%
^{n},\Vert \cdot \Vert _{(2)}\right) $, $\mathbb{S}^{n-1}=\{\alpha \in
\mathbb{R}^{n}\left\vert \Vert \alpha \Vert _{\left( 2\right) }=1\right. {\ }%
\}$ be the unit sphere in $\mathbb{R}^{n}$, $B_{\left( 2\right)
}^{n}=\{\alpha \in \mathbb{R}^{n}\left\vert \Vert \alpha \Vert _{(2)}\leq
1\right. {\ }\}$ be the canonical unit ball in $\mathbb{R}^{n}$ and $\mathrm{%
Vol}_{n}$ be the standard $n$-dimensional volume of subsets in $\mathbb{R}%
^{n}$. Let $V$ be a convex origin-symmetric (i.e. $V=-V$) body in $\mathbb{R}%
^{n}$. Fix a norm $\Vert \cdot \Vert _{V}$ on $\mathbb{R}^{n}$ and denote by
$E$ the Banach space $E=(\mathbb{R}^{n},\Vert \cdot \Vert _{V})$ with the
unit ball $V$. We will use the expectation of the function $\alpha \mapsto
\left\Vert \alpha \right\Vert _{V}$ on $\mathbb{S}^{n-1}$, i.e.%
\[
\mathbf{E}\left[ \left\Vert \cdot \right\Vert _{V}\right] =\int_{\mathbb{S}%
^{n-1}}\left\Vert \alpha \right\Vert _{V}d\mu \left( \alpha \right) ,
\]%
where $d\mu \left( \alpha \right) $ denotes the Haar measure on $\mathbb{S}%
^{n-1}$. The Minkowski sum of two sets $V$ and $W$ in $\mathbb{R}^{n}$ is
defined as%
\[
V+W=\left\{ x+y\left\vert x\in V, y\in W\right. \right\} .
\]%
Let $V\subset \mathbb{R}^{n}$ be an origin-symmetric set. Define its radius $%
\mathrm{rad}\left( V\left\vert W\right. \right) $ in $(\mathbb{R}^{n},\Vert
\cdot \Vert _{W})$ as
\begin{equation}
\mathrm{rad}\left( V\left\vert W\right. \right) =\sup \left\{ \left\Vert
\alpha \right\Vert _{W}\left\vert \alpha \in V\right. \right\} .  \label{8}
\end{equation}

Let $\Phi \left( n\right) =\mathrm{lin}\left\{ \phi _{1},\cdots ,\phi
_{n}\right\} \subset L_{p}$ be a fixed orthonormal system on $\left( \Omega
,\Sigma ,\nu \right) $ and \textrm{J} be the coordinate isomorphism,%
\[
\begin{array}{ccc}
\mathrm{J}:\mathbb{R}^{n} & \longrightarrow & \Phi \left( n\right) \\
\alpha & \longmapsto & t^{\alpha }=\sum_{k=1}^{n}\alpha _{k}\phi _{k}.%
\end{array}%
\]%
The definition $\left\Vert \alpha \right\Vert _{\left( \mathrm{J},p\right)
}:=\left\Vert t^{\alpha }\right\Vert _{p}$ induces a norm on $\mathbb{R}^{n}$%
. Observe that the set
\[
B_{\left( \mathrm{J},p\right) }^{n}=\left\{ \alpha \left\vert \alpha \in
\mathbb{R}^{n}, \left\Vert \alpha \right\Vert _{\left( \mathrm{J}%
,p\right) }\leq 1\right. \right\}
\]%
is a convex and origin-symmetric body in $\mathbb{R}^{n}$ and $B_{\left(
\mathrm{J},2\right) }^{n}=B_{\left( 2\right) }^{n}$.

Various methods of asymptotic geometric analysis have been recently applied
to study the problem of numerical integration on general domains \cite{H-1},
\cite{H-2}, \cite{ku-jfa}. Let $\Lambda _{n}=\mathrm{diag}\left\{ \lambda
_{1},\cdots ,\lambda _{n}\right\} $. The radius $\mathrm{rad}\left( \Lambda
_{n}B_{\left( 2\right) }^{n}\cap L_{s}\left\vert B_{\left( 2\right)
}^{n}\cap L_{s}\right. \right) $, defined by (\ref{8}) of a random
intersetion $\Lambda _{n}B_{\left( 2\right) }^{n}\cap L^{m}$ of
origin-symmetric ellipsoid $\Lambda _{n}B_{\left( 2\right) }^{n}$ with half
exes $\lambda _{1}\geq \cdots \geq \lambda _{n}$ with a random subspace $%
L^{m}\subset \mathbb{R}^{n}$, \textrm{codim}$L^{m}=m$ has been estimated in
\cite{NOVAK1} , \cite{NOVAK1-1}, \cite{NOVAK2}. It was shown that%
\[
\mathrm{rad}\left( \Lambda _{n}B_{\left( 2\right) }^{n}\cap L^{m}\left\vert
B_{\left( 2\right) }^{n}\cap L^{m}\right. \right) \leq Cm^{-\frac{1}{2}%
}\left( \sum_{j\geq \frac{m}{4}}\lambda _{j}^{2}\right) ^{\frac{1}{2}}
\]%
with overwhelming probability, where $C>0$ is an absolute constant.

In Section \ref{SECTIONS} we establish estimates of volumes of sections of
convex origin-symmetric bodies $\mathbf{A}B_{\left( \mathrm{J},p\right)
}^{n} $, where $\mathbf{A}$ is $n\times n$ matrix, $\det \mathbf{A\neq }0$.
These estimates involve expectations $\mathbf{E}\left[ \left\Vert \cdot
\right\Vert _{B_{\left( \mathrm{J},p\right) }^{n}}\right] $ which have been
obtained in \cite{ku-tas} and estimates of volume $\mathrm{Vol}_{n}\left(
B_{\left( \mathrm{J},1\right) }^{n}\right) $ found in \cite{ellipsoid}. We
apply these volume estimates in the proof of Theorem \ref{section-22} and
Theorem \ref{theorem-33}. More precisely, we give lower bounds for the radii
of sections $\mathbf{A}B_{\left( \mathrm{J},p\right) }^{n}\cap L_{r}$ in $%
B_{\left( \mathrm{J},q\right) }^{n}\cap L_{r}$, $\forall L_{r}\subset
\mathbb{R}^{n},1\leq q\leq 2\leq p<\infty $. These estimates are sufficient
to get sharp orders of Gelfand and Kolmogorov widths in Section \ref%
{EXAMPLES} where we consider applications of Theorem \ref{section-22} and
Theorem \ref{theorem-33} in the case $\left( \Omega ,\mathcal{A},\nu \right)
$, where $\Omega =\mathbb{M}^{d}$ is a compact homogeneous Riemannian
manifold and $\nu =\upsilon $ is the Haar measure. Important examples of
such manifolds are given by real and complex Grassmannians, $n$-torus,
Stiefel manifolds, complex spheres and two-point homogeneous spaces, i.e.
real spheres, real, complex and quaternionic projective spaces and the
Cayley elliptic plane ($\mathbb{S}^{d}$, $\mathrm{P}^{d}\left( \mathbb{R}%
\right) $, $\mathrm{P}^{d}\left( \mathbb{C}\right) $, $\mathrm{P}^{d}\left(
\mathbb{H}\right) $, $\mathrm{P}^{16}\left( \mathrm{Cay}\right) $
respectively).

In particular, in the case of two-point homogeneous spaces we show that
\[
\lambda ^{n}\left( W_{p}^{\gamma },L_{1}\right) \asymp n^{-\frac{\gamma }{d}}%
, \gamma >0, p<\infty
\]%
and%
\[
d_{n}\left( W_{\infty }^{\gamma },L_{q}\right) \asymp n^{-\frac{\gamma }{d}}%
, \gamma >0, 1<q\leq 2\leq p\leq \infty .
\]%
These estimates extend the result
\[
d_{n}\left( W_{p}^{\gamma },L_{q}\right) \asymp n^{-\frac{\gamma }{d}}
, \gamma >0, 1<q\leq 2\leq p<\infty
\]%
obtained in \cite{ku-tas}. Also, we show that for Sobolev classes $%
W_{p}^{\gamma }$ on a compact, homogeneous Riemannian manifold we have%
\[
\lambda ^{n}\left( W_{p}^{\gamma },L_{q}\right) \asymp n^{-\frac{\gamma }{d}}%
, \gamma >0, 1<q\leq 2\leq p<\infty .
\]

The method's possibilities are not confined to the statements proved but can
be applied in studying more general problems.

We use several universal constants which enter into the estimates. These
positive constants are denoted by $C_{1}$, $C_{2}$ etc. We did not carefully
distinguish between the different constants, neither did try to get good
estimates for them. Different constants are denoted by the same symbol $C$, $%
C_{1}$, $C_{2}$ etc in different parts of this article. $C_{p}$ denotes the
constant which depends just on $p$. For easy of notation we will write $%
a_{n}\ll b_{n}$ for two sequences if $a_{n}\leq Cb_{n}$ for all $n\in
\mathbb{N}$ and $a_{n}\asymp b_{n}$ if $C_{1}b_{n}\leq a_{n}\leq C_{2}b_{n}$
for all $n\in \mathbb{N}$ and some constants $C_{1}$ and $C_{2}$. Throughout
the text $\left[ a\right] $ denotes the integer part of $a\in \mathbb{R}$.

\section{Volumes of sections of convex origin-symmetric bodies}

\label{SECTIONS}

The next statements give volume estimates which are sufficient for our
applications. Let $V\subset \mathbb{R}^{n}$ and $\mathrm{P}\left(
L_{s}\right) :$ $\mathbb{R}^{n}\rightarrow L_{s}$ be the operator of
orthogonal projection onto $L_{s}\subset \mathbb{R}^{n}$. Let $\mathbf{O}%
\left( n\right) $ be the orthogonal group in $\mathbb{R}^{n}$. Then the
Grassmannian $\mathbb{G}\left( s,n\right) =\mathbf{O}\left( n\right) /\left(
\mathbf{O}\left( s\right) \times \mathbf{O}\left( n-s\right) \right) $ is
the space of $s$ dimensional subspaces $L_{s}$ in $\mathbb{R}^{n}$. Define

\[
\begin{array}{ccc}
\Upsilon \left( L_{s},V\right) :\mathbb{G}\left( s,n\right) & \longrightarrow
& \mathbb{R} \\
L_{s} & \longmapsto & \mathrm{Vol}_{s}\left( \mathrm{P}\left( L_{s}\right)
V\right) .%
\end{array}%
\]

{\bf Lemma 1}
\label{el-section}
\begin{em}
Let $\mathbf{A}$ be $n\times n$ matrix, $\det \mathbf{%
A\neq }0$. Then%
\[
\max \left\{ \Upsilon \left( L_{s},\mathbf{A}B_{\left( 2\right) }^{n}\right)
\left\vert L_{s}\in \mathbb{G}\left( s,n\right) \right. \right\} \leq
3^{n}\varrho _{n}^{s-n}\det \left( \mathbf{A}\right) \mathrm{Vol}_{s}\left(
B_{\left( 2\right) }^{s}\right) ,
\]%
where%
\[
\varrho _{n}=\left( \rho \left( \left( \mathbf{A}^{-1}\right) ^{\ast }%
\mathbf{A}^{-1}\right) \right) ^{-\frac{1}{2}}
\]%
and
\[
\rho \left( \left( \mathbf{A}^{-1}\right) ^{\ast }\mathbf{A}^{-1}\right)
=\max \left\{ \left\vert \sigma _{k}\right\vert , 1\leq k\leq
n\right\}
\]%
is the spectral radius of $\left( \mathbf{A}^{-1}\right) ^{\ast }\mathbf{A}%
^{-1}$ with the eigenvalues $\sigma _{k}$, $1\leq k\leq n$.
\end{em}

{\bf Proof}
Recall that $\left\{ x_{j}, 1\leq j\leq n\left( \delta \right)
\right\} \subset V$ is called a $\delta $-net of a compact set $V\subset
\mathbb{R}^{n}$ in $(\mathbb{R}^{n},\Vert \cdot \Vert _{W})$ if for each $%
x\in V$ there is one $x_{j}\in V$ such that $\left\Vert x-x_{j}\right\Vert
_{W}\leq \delta $ or%
\[
V\subset \cup _{j=1}^{n\left( \delta \right) }\left( x_{j}+\delta W\right) .
\]%
Also, a set of points $\left\{ y_{j}, 1\leq j\leq m\left( \delta
\right) \right\} \subset V$ is called $\delta $-distinguishable in $\left(
\mathbb{R}^{n},\left\Vert \cdot \right\Vert _{W}\right) $ if
\[
\left\Vert y_{k}-y_{j}\right\Vert _{W}\geq \delta , \forall 1\leq
k\neq j\leq m\left( \delta \right) .
\]%
Let $\left\{ x_{j}, 1\leq j\leq n\left( \delta \right) \right\} $ be
a $\delta $-net of minimal cardinality $n\left( \delta \right) $ and $%
\left\{ y_{j}, 1\leq j\leq m\left( \delta \right) \right\} $ be a $%
\delta $-distinguishable net of maximal cardinality $m\left( \delta \right) $%
. It is well-known that
\begin{equation}
m\left( 2\delta \right) \leq n\left( \delta \right) \leq m\left( \delta
\right) ,  \label{mmm}
\end{equation}%
(see e.g. \cite{mak}). Let $V=\mathbf{A}B_{\left( 2\right) }^{n}$ and $%
W=B_{\left( 2\right) }^{n}$ then the sets
\begin{equation}
A_{j}=y_{j}+\frac{\delta }{2}B_{\left( 2\right) }^{n}, 1\leq j\leq
m\left( \delta \right)  \label{delta1}
\end{equation}%
are disjoint. It is well-known that%
\[
\left\Vert \mathbf{A}^{-1}\left\vert l_{2}^{n}\rightarrow l_{2}^{n}\right.
\right\Vert =\left( \rho _{n}\left( \left( \mathbf{A}^{-1}\right) ^{\ast }%
\mathbf{A}^{-1}\right) \right) ^{\frac{1}{2}}.
\]
Hence%
\[
\left\Vert \mathbf{A}^{-1}x\right\Vert _{\left( 2\right) }\leq \left( \rho
_{n}\left( \left( \mathbf{A}^{-1}\right) ^{\ast }\mathbf{A}^{-1}\right)
\right) ^{\frac{1}{2}}\left\Vert x\right\Vert _{\left( 2\right) },
\]%
or $\varrho _{n}B_{\left( 2\right) }^{n}\subset \mathbf{A}B_{\left( 2\right)
}^{n}$, where $\varrho _{n}=\left( \rho _{n}\left( \left( \mathbf{A}%
^{-1}\right) ^{\ast }\mathbf{A}^{-1}\right) \right) ^{-\frac{1}{2}}$. Let $%
\delta =\varrho _{n}$ in (\ref{delta1}). Then
\[
A_{j}\subset \mathbf{A}B_{\left( 2\right) }^{n}+\frac{\varrho _{n}}{2}%
B_{\left( 2\right) }^{n}
\]%
\[
\subset \mathbf{A}B_{\left( 2\right) }^{n}+\frac{1}{2}\mathbf{A}B_{\left(
2\right) }^{n}=\frac{3}{2}\mathbf{A}B_{\left( 2\right) }^{n}, \forall
1\leq j\leq m\left( \varrho _{n}\right) .
\]%
Consequently,
\[
\cup _{j=1}^{m\left( \varrho _{n}\right) }A_{j}\subset \frac{3}{2}\mathbf{A}%
B_{\left( 2\right) }^{n}.
\]%
By comparing volumes we get%
\[
m\left( \varrho _{n}\right) \left( \frac{\varrho _{n}}{2}\right) ^{n}\mathrm{%
Vol}_{n}\left( B_{\left( 2\right) }^{n}\right) \leq \left( \frac{3}{2}%
\right) ^{n}\det \left( \mathbf{A}\right) \mathrm{Vol}_{n}\left( B_{\left(
2\right) }^{n}\right) ,
\]%
or%
\begin{equation}
m\left( \varrho _{n}\right) \leq 3^{n}\varrho _{n}^{-n}\det \left( \mathbf{A}%
\right) .  \label{nnn}
\end{equation}%
From (\ref{mmm}) and (\ref{nnn}) we find%
\begin{equation}
n\left( \varrho _{n}\right) \leq m\left( \varrho _{n}\right) \leq
3^{n}\varrho _{n}^{-n}\det \left( \mathbf{A}\right) .  \label{ineq}
\end{equation}%
This means that
\begin{equation}
\mathbf{A}B_{\left( 2\right) }^{n}\subset \cup _{j=1}^{n\left( \varrho
_{n}\right) }\left( x_{j}+\varrho _{n}B_{\left( 2\right) }^{n}\right) .
\label{incl}
\end{equation}%
Let $L_{s}$ be an arbitrary subspace in $\mathbb{R}^{n}$, $\dim L_{s}=s$.
Then applying (\ref{incl}) and (\ref{ineq}) we obtain
\[
\mathrm{Vol}_{s}\left( \mathrm{P}\left( L_{s}\right) \mathbf{A}B_{\left(
2\right) }^{n}\right) \leq n\left( \varrho _{n}\right) \varrho _{n}^{s}%
\mathrm{Vol}_{s}\left( \mathrm{P}\left( L_{s}\right) B_{\left( 2\right)
}^{n}\right)
\]%
\[
=n\left( \varrho _{n}\right) \varrho _{n}^{s}\mathrm{Vol}_{s}\left(
B_{\left( 2\right) }^{s}\right)
\]%
\[
\leq 3^{n}\varrho _{n}^{s-n}\det \left( \mathbf{A}\right) \mathrm{Vol}%
_{s}\left( B_{\left( 2\right) }^{s}\right) .
\]
$\Box$

{Lemma 2}
\label{lemma-2}
\begin{em}
(a) Let $\Phi \left( n\right) =\mathrm{lin}\left\{ \phi _{1},\cdots ,\phi
_{n}\right\} \subset L_{p}$, $1<p\leq 2$, $n\in \mathbb{N}$ be a fixed
orthonormal system on $\left( \Omega ,\Sigma ,\nu \right) $ then%
\[
\max \left\{ \Upsilon \left( L_{s},B_{\left( \mathrm{J},p\right)
}^{n}\right) \left\vert L_{s}\in \mathbb{G}\left( s,n\right) \right.
\right\}
\]%
\begin{equation}
\leq \left( \frac{5}{2}\right) ^{n}\left( \mathbf{E}\left[ \left\Vert \cdot
\right\Vert _{B_{\left( \mathrm{J},p^{^{\prime }}\right) }^{n}}\right]
\right) ^{n}\mathrm{Vol}_{s}\left( B_{\left( 2\right) }^{s}\right).
\label{pppp}
\end{equation}%
(b) Let $\left\Vert \phi _{k}\right\Vert _{\infty }\leq M$, $1\leq k\leq n$
then%
\[
\max \left\{ \Upsilon \left( L_{s},B_{\left( \mathrm{J},1\right)
}^{n}\right) \left\vert L_{s}\in \mathbb{G}\left( s,n\right) \right.
\right\} \leq C^{n}\mathrm{Vol}_{s}\left( B_{\left( 2\right) }^{s}\right) .
\]
\end{em}

{Proof}
We prove (b) first. Let $\left\{ x_{j}, 1\leq j\leq n\left( \frac{1}{2%
}\right) \right\} $ be a $\frac{1}{2}$-net of minimal cardinality $n\left(
\frac{1}{2}\right) $ for $B_{\left( \mathrm{J},1\right) }^{n}$ in $(\mathbb{R%
}^{n},\Vert \cdot \Vert _{B_{\left( 2\right) }^{n}})$ and $\left\{ y_{j}, 1\leq j\leq m\left( \frac{1}{2}\right) \right\} $ be a $\frac{1}{2}$%
-distinguishable net for $B_{\left( \mathrm{J},1\right) }^{n}$ in $(\mathbb{R%
}^{n},\Vert \cdot \Vert _{B_{\left( 2\right) }^{n}})$ of maximal cardinality
$m\left( \frac{1}{2}\right) $. Then the sets
\[
A_{j}=y_{j}+\frac{1}{4}B_{\left( 2\right) }^{n}, 1\leq j\leq m\left(
\frac{1}{2}\right)
\]%
are disjoint and
\[
A_{j}\subset B_{\left( \mathrm{J},1\right) }^{n}+\frac{1}{4}B_{\left(
2\right) }^{n}\subset B_{\left( \mathrm{J},1\right) }^{n}+\frac{1}{4}%
B_{\left( \mathrm{J},1\right) }^{n}
\]%
\[
=\frac{5}{4}B_{\left( \mathrm{J},1\right) }^{n}, \forall 1\leq j\leq
m\left( \frac{1}{2}\right)
\]%
since $B_{\left( 2\right) }^{n}\subset B_{\left( \mathrm{J},1\right) }^{n}$.
Consequently,
\[
\cup _{j=1}^{m\left( \frac{1}{2}\right) }A_{j}\subset \frac{5}{4}B_{\left(
\mathrm{J},1\right) }^{n}.
\]%
By comparing volumes we get%
\[
m\left( \frac{1}{2}\right) \left( \frac{1}{4}\right) ^{n}\mathrm{Vol}%
_{n}\left( B_{\left( 2\right) }^{n}\right) \leq \left( \frac{5}{4}\right)
^{n}\mathrm{Vol}_{n}\left( B_{\left( \mathrm{J},1\right) }^{n}\right) ,
\]%
or
\begin{equation}
m\left( \frac{1}{2}\right) \leq 5^{n}\frac{\mathrm{Vol}_{n}\left( B_{\left(
\mathrm{J},1\right) }^{n}\right) }{\mathrm{Vol}_{n}\left( B_{\left( 2\right)
}^{n}\right) }.  \label{vol-v}
\end{equation}%
It is shown in \cite{ellipsoid} Lemma 3, (26) (see also \cite{klw}) that
\begin{equation}
\mathrm{Vol}_{n}\left( B_{\left( \mathrm{J},1\right) }^{n}\right) \leq C^{n}%
\mathrm{Vol}_{n}\left( B_{\left( 2\right) }^{n}\right) .  \label{b1}
\end{equation}%
Hence, by (\ref{mmm}), (\ref{vol-v}) and (\ref{b1}) we get
\begin{equation}
n\left( \frac{1}{2}\right) \leq m\left( \frac{1}{2}\right) \leq C^{n}.
\label{ccnn}
\end{equation}%
This means that there are such $x_{j},$ $1\leq j\leq n\left( \frac{1}{2}%
\right) $ that
\begin{equation}
B_{\left( \mathrm{J},1\right) }^{n}\subset \cup _{j=1}^{n\left( \frac{1}{2}%
\right) }\left( x_{j}+\frac{1}{2}B_{\left( 2\right) }^{n}\right) .
\label{cc11}
\end{equation}%
Let $L_{s}$ be an arbitrary subspace in $\mathbb{R}^{n}$, $\dim L_{s}=s$,
then applying (\ref{ccnn}) and (\ref{cc11}) we find
\[
\mathrm{Vol}_{s}\left( \mathrm{P}\left( L_{s}\right) B_{\left( \mathrm{J}%
,1\right) }^{n}\right) \leq 2^{-n}n\left( \frac{1}{2}\right) \mathrm{Vol}%
_{s}\left( \mathrm{P}\left( L_{s}\right) B_{\left( 2\right) }^{n}\right)
\]%
\[
=2^{-n}n\left( \frac{1}{2}\right) \mathrm{Vol}_{s}\left( B_{\left( 2\right)
}^{s}\right) \leq C^{n}\mathrm{Vol}_{s}\left( B_{\left( 2\right)
}^{s}\right) .
\]%
Let us show (a) now. Let $V$ be an origin-symmetric body in $\mathbb{R}^{n}$
which is the unit ball in $E=(\mathbb{R}^{n},\Vert \cdot \Vert _{V})$. A
direct calculation (see e.g. \cite{pisier}, p.xi, \cite{ku-tas}) shows that
\begin{equation}
\frac{\mathrm{Vol}_{n}\left( V\right) }{\mathrm{Vol}_{n}\left( B_{\left(
2\right) }^{n}\right) }=\int_{\mathbb{S}^{n-1}}\left\Vert \alpha \right\Vert
_{V}^{-n}d\mu \left( \alpha \right)  \label{111}
\end{equation}%
and by convexity we get%
\begin{equation}
\int_{\mathbb{S}^{n-1}}\left\Vert \alpha \right\Vert _{V}^{-n}d\mu \left(
\alpha \right) \geq \left( \int_{\mathbb{S}^{n-1}}\left\Vert \alpha
\right\Vert _{V}d\mu \left( \alpha \right) \right) ^{-n}=\mathbf{E}\left[
\left\Vert \cdot \right\Vert _{V}\right] ^{-n}.  \label{222}
\end{equation}%
From (\ref{111}) and (\ref{222}) we find%
\begin{equation}
\mathrm{Vol}_{n}\left( V\right) \geq \left( \mathbf{E}\left[ \left\Vert
\cdot \right\Vert _{V}\right] \right) ^{-n}\mathrm{Vol}_{n}\left( B_{\left(
2\right) }^{n}\right) .  \label{333}
\end{equation}%
Comparing (\ref{333}) with Santalo inequality \cite{santalo},
\[
\frac{\mathrm{Vol}_{n}\left( V\right) \mathrm{Vol}_{n}\left( V^{o}\right) }{%
\left( \mathrm{Vol}_{n}\left( B_{\left( 2\right) }^{n}\right) \right) ^{2}}%
\leq 1,
\]%
which is valid for any convex and origin-symmetric body $V\subset \mathbb{R}%
^{n}$, we get
\begin{equation}
\mathrm{Vol}_{n}\left( V\right) \leq \left( \mathbf{E}\left[ \left\Vert
\cdot \right\Vert _{V^{o}}\right] \right) ^{n}\mathrm{Vol}_{n}\left(
B_{\left( 2\right) }^{n}\right) .  \label{vol-v1}
\end{equation}%
In particular, let $V=B_{\left( \mathrm{J},p\right) }^{n}$, $1<p\leq 2$ in (%
\ref{vol-v1}) then%
\begin{equation}
\mathrm{Vol}_{n}\left( B_{\left( \mathrm{J},p\right) }^{n}\right) \leq
\left( \mathbf{E}\left[ \left\Vert \cdot \right\Vert _{\left( B_{\left(
\mathrm{J},p\right) }^{n}\right) ^{o}}\right] \right) ^{n}\mathrm{Vol}%
_{n}\left( B_{\left( 2\right) }^{n}\right) .  \label{vol-p}
\end{equation}%
Since
\[
\mathbf{E}\left[ \left\Vert \cdot \right\Vert _{\left( B_{\left( \mathrm{J}%
,p\right) }^{n}\right) ^{o}}\right] \leq \mathbf{E}\left[ \left\Vert \cdot
\right\Vert _{B_{\left( \mathrm{J},p^{^{\prime }}\right) }^{n}}\right] 
, \frac{1}{p}+\frac{1}{p^{^{\prime }}}=1
\]%
and $B_{\left( 2\right) }^{n}\subset B_{\left( \mathrm{J},p\right) }^{n}$if $%
1<p\leq 2$, then form (\ref{vol-p}) and we obtain
\[
m\left( \frac{1}{2}\right) \leq 5^{n}\frac{\mathrm{Vol}_{n}\left( B_{\left(
\mathrm{J},p\right) }^{n}\right) }{\mathrm{Vol}_{n}\left( B_{\left( 2\right)
}^{n}\right) }
\]%
\[
\leq 5^{n}\left( \mathbf{E}\left[ \left\Vert \cdot \right\Vert _{\left(
B_{\left( \mathrm{J},p\right) }^{n}\right) ^{o}}\right] \right) ^{n}\leq
5^{n}\left( \mathbf{E}\left[ \left\Vert \cdot \right\Vert _{B_{\left(
\mathrm{J},p^{^{\prime }}\right) }^{n}}\right] \right) ^{n}.
\]%
Hence by (\ref{mmm}),%
\[
n\left( \frac{1}{2}\right) \leq m\left( \frac{1}{2}\right) \leq 5^{n}\left(
\mathbf{E}\left[ \left\Vert \cdot \right\Vert _{B_{\left( \mathrm{J}%
,p^{^{\prime }}\right) }^{n}}\right] \right) ^{n}.
\]%
Since%
\[
B_{\left( \mathrm{J},p\right) }^{n}\subset \cup _{j=1}^{n\left( \frac{1}{2}%
\right) }\left( x_{j}+\frac{1}{2}B_{\left( 2\right) }^{n}\right)
\]%
then%
\[
\mathrm{Vol}_{s}\left( \mathrm{P}\left( L_{s}\right) B_{\left( \mathrm{J}%
,p\right) }^{n}\right) \leq 2^{-n}n\left( \frac{1}{2}\right) \mathrm{Vol}%
_{s}\left( B_{\left( 2\right) }^{s}\right)
\]%
\[
\leq \left( \frac{5}{2}\right) ^{n}\left( \mathbf{E}\left[ \left\Vert \cdot
\right\Vert _{B_{\left( \mathrm{J},p^{^{\prime }}\right) }^{n}}\right]
\right) ^{n}\mathrm{Vol}_{s}\left( B_{\left( 2\right) }^{s}\right) .
\]
$\Box$

We will need the following definition from \cite{ellipsoid}.

{Definition 1}
\begin{em}
We say that an orthonormal system $\left\{ \phi _{k}, k\in \mathbb{N}%
\right\} $ on a probability space $\left( \Omega ,\Sigma ,\nu \right) $ is
essentially bounded if for any $0<C<1$ and $n\in \mathbb{N}$ there exists a
set of orthonormal functions $\left\{ \xi _{k}, 1\leq k\leq m\right\}
$, $m\geq Cn$ such that%
\[
\Xi \left( m\right) =\mathrm{lin}\left\{ \xi _{k}, 1\leq k\leq
m\right\} \subset \Phi \left( n\right) =\mathrm{lin}\left\{ \phi _{k}, 1\leq k\leq n\right\}
\]%
and
\[
\max \left\{ \left\Vert \xi _{k}\right\Vert _{\infty }\left\vert 1\leq k\leq
m\right. \right\} \leq M,
\]%
where $M>0$ dependes just on $C$.
\end{em}

The next statement is a consequence of Lemma \ref{el-section} and Lemma \ref%
{lemma-2}.

{Theorem 1}
\begin{em}
\label{section-22} Let $\left\{ \phi _{k}, k\in \mathbb{N}\right\} $
be an essentially bounded orthonormal system, $\mathbf{A}$ be $n\times n$
matrix, $\det \mathbf{A\neq }0$ and $r\geq \frac{2}{3}n$ then
\[
\min \left\{ \mathrm{rad}\left( \mathbf{A}B_{\left( \mathrm{J},p\right)
}^{n}\cap L_{r}\left\vert B_{\left( \mathrm{J},1\right) }^{n}\cap
L_{r}\right. \right) \left\vert L_{r}\in \mathbb{G}\left( r,n\right) \right.
\right\}
\]%
\[
\geq C\varrho _{n}\mathbf{E}\left[ \left\Vert \cdot \right\Vert _{B_{\left(
\mathrm{J},p\right) }^{n}}\right] ^{-\frac{3}{2}}, p\geq 2,
\]%
where
\[
\varrho _{n}=\left( \rho \left( \left( \mathbf{A}^{-1}\right) ^{\ast }%
\mathbf{A}^{-1}\right) \right) ^{-\frac{1}{2}}
\]%
and
\[
\rho \left( \left( \mathbf{A}^{-1}\right) ^{\ast }\mathbf{A}^{-1}\right)
=\max \left\{ \left\vert \sigma _{k}\right\vert , 1\leq k\leq
n\right\}
\]%
is the spectral radius of $\left( \mathbf{A}^{-1}\right) ^{\ast }\mathbf{A}%
^{-1}$ with the eigenvalues $\sigma _{k}$, $1\leq k\leq n$.
\end{em}

{\bf Proof}
Let $L_{s}\subset \mathbb{R}^{n}$, $1\leq s\leq n-1$ be an arbitrary
subspace, $\dim L_{s}=s$, $V\subset \mathbb{R}^{n}$ be an origin-symmetric
convex body and $L_{s}^{\bot }$ denotes the orthogonal complement of $L_{s}$
in $\mathbb{R}^{n}$. Then%
\begin{equation}
\mathrm{Vol}_{n}\left( V\right) =\int_{V}dz=\int_{\mathrm{P}\left(
L_{s}^{\bot }\right) V}\mathrm{Vol}_{s}\left( V\cap \left( z+L_{s}\right)
\right) dz.  \label{00}
\end{equation}%
Hence, involving standard arguments connected with the Brunn-Minkowski
theorem (see e.g. \cite{bourgain}, \cite{ku-jfa}, \cite{ku-tas}) we get
\begin{equation}
\mathrm{Vol}_{s}\left( V\cap \left( z+L_{s}\right) \right) \leq \mathrm{Vol}%
_{s}\left( V\cap L_{s}\right) , \forall z\in \mathrm{P}\left(
L_{s}^{\bot }\right) V.  \label{000}
\end{equation}%
Let $V=\mathbf{A}B_{\left( \mathrm{J},p\right) }^{n}$, $p\geq 2$. Since $%
B_{\left( \mathrm{J},p\right) }^{n}\subset B_{\left( 2\right) }^{n}$ then $%
\mathbf{A}B_{\left( \mathrm{J},p\right) }^{n}\subset \mathbf{A}B_{\left(
2\right) }^{n}$. Hence by (\ref{00}), (\ref{000}) and Lemma \ref{el-section}%
,
\[
\mathrm{Vol}_{n}\left( \mathbf{A}B_{\left( \mathrm{J},p\right) }^{n}\right)
\leq \mathrm{Vol}_{s}\left( \mathbf{A}B_{\left( \mathrm{J},p\right)
}^{n}\cap L_{s}\right) \mathrm{Vol}_{n-s}\left( \mathrm{P}\left( L_{s}^{\bot
}\right) \mathbf{A}B_{\left( \mathrm{J},p\right) }^{n}\right)
\]%
\[
\leq \mathrm{Vol}_{s}\left( \mathbf{A}B_{\left( \mathrm{J},p\right)
}^{n}\cap L_{s}\right) \mathrm{Vol}_{n-s}\left( \mathrm{P}\left( L_{s}^{\bot
}\right) \mathbf{A}B_{\left( 2\right) }^{n}\right)
\]%
\begin{equation}
\leq \mathrm{Vol}_{s}\left( \mathbf{A}B_{\left( \mathrm{J},p\right)
}^{n}\cap L_{s}\right) 3^{n}\varrho _{n}^{-s}\det \left( \mathbf{A}\right)
\mathrm{Vol}_{n-s}\left( B_{\left( 2\right) }^{n-s}\right) .  \label{aaa}
\end{equation}%
Applying (\ref{333}) we find
\[
\mathrm{Vol}_{n}\left( \mathbf{A}B_{\left( \mathrm{J},p\right) }^{n}\right)
=\det \left( \mathbf{A}\right) \mathrm{Vol}_{n}\left( B_{\left( \mathrm{J}%
,p\right) }^{n}\right)
\]%
\begin{equation}
\geq \det \left( \mathbf{A}\right) \mathbf{E}\left[ \left\Vert \cdot
\right\Vert _{B_{\left( \mathrm{J},p\right) }^{n}}\right] ^{-n}\mathrm{Vol}%
_{n}\left( B_{\left( 2\right) }^{n}\right) .  \label{bbb}
\end{equation}%
Comparing (\ref{aaa}) and (\ref{bbb}) we obtain%
\begin{equation}
\mathrm{Vol}_{s}\left( \mathbf{A}B_{\left( \mathrm{J},p\right) }^{n}\cap
L_{s}\right) \geq \frac{\varrho _{n}^{s}\mathrm{Vol}_{n}\left( B_{\left(
2\right) }^{n}\right) }{3^{n}\left( \mathbf{E}\left[ \left\Vert \cdot
\right\Vert _{B_{\left( \mathrm{J},p\right) }^{n}}\right] \right) ^{n}%
\mathrm{Vol}_{n-s}\left( B_{\left( 2\right) }^{n-s}\right) }.  \label{ccc}
\end{equation}%
Since the orthonormal system $\left\{ \phi _{k}, k\in \mathbb{N}%
\right\} $ is essentially bounded then for any $n\in \mathbb{N}$ there
exists an orthonormal system $\left\{ \xi _{k}, 1\leq k\leq m\right\}
$, $m\geq Cn$ such that $\left\Vert \xi _{k}\right\Vert _{\infty }\leq M$, $%
1\leq k\leq m$, where $M$ depends just on $C$. Let $\Xi \left( m\right) =%
\mathrm{lin}\left\{ \xi _{k}, 1\leq k\leq m\right\} $ and $\mathrm{K}$
be the coordinate isomorphism,%
\[
\begin{array}{ccc}
\mathrm{K}:\mathbb{R}^{m} & \longrightarrow & \Xi \left( m\right) \\
\alpha & \longmapsto & t^{\alpha }=\sum_{k=1}^{m}\alpha _{k}\xi _{k}.%
\end{array}%
\]%
The definition $\left\Vert \alpha \right\Vert _{\left( \mathrm{K},p\right)
}:=\left\Vert t^{\alpha }\right\Vert _{p}$ induces a norm on $\mathbb{R}^{m}$%
. Clearly, the set
\[
B_{\left( \mathrm{K},p\right) }^{m}=\left\{ \alpha \left\vert \alpha \in
\mathbb{R}^{m}, \left\Vert \alpha \right\Vert _{\left( \mathrm{K}%
,p\right) }\leq 1\right. \right\}
\]%
is a convex and origin-symmetric body in $\mathbb{R}^{m}$ and $B_{\left(
\mathrm{K},2\right) }^{m}=B_{\left( 2\right) }^{m}$. Let $L_{r}$ be
arbitrary subspace in $\mathbb{R}^{n}$ such that $\dim L_{r}=r\geq \left[
\frac{2n}{3}\right] $. Pick $\mathrm{J}^{-1}\Xi \left( m\right) \subset
\mathbb{R}^{n}$, $\dim \mathrm{J}^{-1}\Xi \left( m\right) =m\geq \left[
\frac{2n}{3}\right] $. Clearly,%
\[
s:=\dim \left( L_{r}\cap \mathrm{J}^{-1}\Xi \left( m\right) \right) \geq
r+m-n\geq \frac{n}{3}-2.
\]%
By Lemma \ref{lemma-2} (b)
\[
\mathrm{Vol}_{s}\left( B_{\mathrm{K},1}^{m}\cap L_{s}\right) \leq C^{m}%
\mathrm{Vol}_{s}\left( B_{\left( 2\right) }^{s}\right)
\]%
Consequently, by (\ref{ccc}) we conclude that
\[
\mathrm{rad}\left( \mathbf{A}B_{\left( \mathrm{J},p\right) }^{n}\cap
L_{s}\left\vert B_{\left( \mathrm{J},1\right) }^{n}\cap L_{s}\right. \right)
\geq \left( \frac{\mathrm{Vol}_{s}\left( \mathbf{A}B_{\left( \mathrm{J}%
,p\right) }^{n}\cap L_{s}\right) }{\mathrm{Vol}_{s}\left( B_{\mathrm{K}%
,1}^{m}\cap L_{s}\right) }\right) ^{\frac{1}{s}}
\]%
\[
\geq C^{-\frac{m}{s}}3^{-\frac{n}{s}}\varrho _{n}\left( \mathbf{E}\left[
\left\Vert \cdot \right\Vert _{B_{\left( \mathrm{J},p\right) }^{n}}\right]
\right) ^{-\frac{n}{s}}\left( \frac{\mathrm{Vol}_{n}\left( B_{\left(
2\right) }^{n}\right) }{\mathrm{Vol}_{s}\left( B_{\left( 2\right)
}^{s}\right) ^{n}\mathrm{Vol}_{n-s}\left( B_{\left( 2\right) }^{n-s}\right) }%
\right) ^{\frac{1}{s}}
\]%
\begin{equation}
\geq C\varrho _{n}\mathbf{E}\left[ \left\Vert \cdot \right\Vert _{B_{\left(
\mathrm{J},p\right) }^{n}}\right] ^{-\frac{n}{s}},  \label{ddd}
\end{equation}%
since by Lemma 1, \cite{ku-tas},
\begin{equation}
C_{1}\leq \left( \frac{\mathrm{Vol}_{n}\left( B_{\left( 2\right)
}^{n}\right) }{\mathrm{Vol}_{s}\left( B_{\left( 2\right) }^{s}\right)
\mathrm{Vol}_{n-s}\left( B_{\left( 2\right) }^{n-s}\right) }\right) ^{\frac{1%
}{s}}\leq C_{2}  \label{zzz}
\end{equation}%
if $C_{3}n\leq s\leq C_{4}n$, where $0<C_{3}<C_{4}<1$.
$\Box$

{Theorem 2}
\label{theorem-33} 
\begin{em}
Let $C_{1}n\leq s\leq C_{2}n$, where $0<C_{1}<C_{2}<1$
and $1<q\leq 2\leq p<\infty $ then%
\[
\mathrm{rad}\left( \mathbf{A}B_{\left( \mathrm{J},p\right) }^{n}\cap
L_{s}\left\vert B_{\left( \mathrm{J},q\right) }^{n}\cap L_{s}\right. \right)
\]%
\[
\geq C\varrho _{n}\left( \mathbf{E}\left[ \left\Vert \cdot \right\Vert
_{B_{\left( \mathrm{J},q^{^{\prime }}\right) }^{n}}\right] \mathbf{E}\left[
\left\Vert \cdot \right\Vert _{B_{\left( \mathrm{J},p\right) }^{n}}\right]
\right) ^{-\frac{n}{s}},
\]%
where $\frac{1}{q}+\frac{1}{q^{^{\prime }}}=1$,%
\[
\varrho _{n}=\left( \rho \left( \left( \mathbf{A}^{-1}\right) ^{\ast }%
\mathbf{A}^{-1}\right) \right) ^{-\frac{1}{2}},
\]%
and
\[
\rho \left( \left( \mathbf{A}^{-1}\right) ^{\ast }\mathbf{A}^{-1}\right)
=\max \left\{ \left\vert \sigma _{k}\right\vert , 1\leq k\leq
n\right\}
\]%
is the spectral radius of $\left( \mathbf{A}^{-1}\right) ^{\ast }\mathbf{A}%
^{-1}$ with the eigenvalues $\sigma _{k}$, $1\leq k\leq n$.
\end{em}

{\bf Proof}
Comparing (\ref{pppp}), (\ref{ccc}) and (\ref{zzz}) we get%
\[
\mathrm{rad}\left( \mathbf{A}B_{\left( \mathrm{J},p\right) }^{n}\cap
L_{s}\left\vert B_{\left( \mathrm{J},q\right) }^{n}\cap L_{s}\right. \right)
\geq \left( \frac{\mathrm{Vol}_{s}\left( \mathbf{A}B_{\left( \mathrm{J}%
,p\right) }^{n}\cap L_{s}\right) }{\mathrm{Vol}_{s}\left( B_{\mathrm{J}%
,q}^{n}\cap L_{s}\right) }\right) ^{\frac{1}{s}}
\]%
\[
\geq C\varrho _{n}\left( \mathbf{E}\left[ \left\Vert \cdot \right\Vert
_{B_{\left( \mathrm{J},q^{^{\prime }}\right) }^{n}}\right] \mathbf{E}\left[
\left\Vert \cdot \right\Vert _{B_{\left( \mathrm{J},p\right) }^{n}}\right]
\right) ^{-\frac{n}{s}}.
\]
$\Box$

\section{Examples and applications}

\label{EXAMPLES}

In this section we consider applications of Theorem \ref{section-22} and
Theorem \ref{theorem-33}. Consider a compact, homogeneous $d$-dimensional
Riemannian manifold $\mathbb{M}^{d}$. Let $g$ its metric tensor, $\upsilon $
its normalized volume element and $\Delta $ its Laplace--Beltrami operator.
In local coordinates $x_{l}$, $1\leq l\leq d$,%
\[
\Delta =-(\overline{g})^{-\frac{1}{2}}\sum_{k}\frac{\partial }{\partial x_{k}%
}\left( g^{jk}\right) (\overline{g})^{\frac{1}{2}}\frac{\partial }{\partial
x_{j}},
\]%
where $g_{jk}=g\left( \frac{\partial }{\partial x_{j}},\frac{\partial }{%
\partial x_{k}}\right) $, $\overline{g}=|\det (g_{jk})|$ , and $\left(
g^{jk}\right) =(g_{jk})^{-1}$. It is well-known that $\Delta $ is an
elliptic, self adjoint, invariant under isometry, second order operator. The
eigenvalues $\theta _{k}$, $k\geq 0$ of $\Delta $ are discrete, nonnegative
and form an increasing sequence with $+\infty $ the only accumulation point.
The corresponding eigenspaces $\mathrm{H}_{k}$, $k\geq 0$ are finite
dimensional, orthogonal with respect to the scalar product
\[
\left[ f,g\right] =\int_{\mathbb{M}^{d}}fgd\upsilon
\]%
and
\[
L_{2}=L_{2}\left( \mathbb{M}^{d},\right) =\oplus _{k=0}^{\infty }\mathrm{H}%
_{k}.
\]%
Fix an orthonormal basis $\left\{ \mathrm{Y}_{m}^{k}\right\} _{m=1}^{d_{k}}$%
, where $d_{k}=\dim \mathrm{H}_{k}$ of $\mathrm{H}_{k}$. A wide range of
sets of smooth functions can be introduced using multiplier operators. Let $%
\varphi \in L_{p}$, $1\leq p\leq \infty $ with the formal Fourier series
\[
\varphi \sim c_{0}+\sum_{k=1}^{\infty }\sum_{m=1}^{d_{k}}c_{k,m}\left(
\varphi \right) \mathrm{Y}_{m}^{k}, c_{k,m}\left( \varphi \right)
=\int_{\mathbb{M}^{d}}\varphi \mathrm{Y}_{m}^{k}d\upsilon
\]%
and $\lambda :\left( 0,\infty \right) \mapsto \mathbb{R}$ be a continuous
function. If for any $\varphi \in L_{p}$ there is a function $f=\Lambda
\varphi \in L_{q}$ such that
\[
f\sim \sum_{k=0}^{\infty }\lambda \left( \theta _{k}\right)
\sum_{m=1}^{\infty }c_{k,m}\left( \varphi \right) \mathrm{Y}_{m}^{k}
\]%
then we say that the multiplier $\Lambda =\left\{ \lambda \left( \theta
_{k}\right) , k\geq 0\right\} $ is of $\left( p,q\right) $-type. In
particular, let $\lambda \left( t\right) =t^{-\frac{\gamma }{2}}$, $\gamma
>0 $ then the $\gamma $-th fractional integral $I_{\gamma }\varphi $ is
defined as%
\begin{equation}
I_{\gamma }\varphi \sim \sum_{k=1}^{\infty }\theta _{k}^{-\frac{\gamma }{2}%
}\sum_{m=1}^{\infty }c_{k,m}\left( \varphi \right) \mathrm{Y}_{m}^{k}.
\label{int-fr}
\end{equation}%
The standard Sobolev classes $W_{p}^{\gamma }$ are defined as sets of
functions with formal Fourier expansions (\ref{int-fr}), where
\[
\left\Vert \varphi \right\Vert _{p}\leq 1, \int_{\mathbb{M}%
^{d}}\varphi d\upsilon =0.
\]%
It is known that (see Lemma 2.2, \cite{K-L-T})%
\begin{equation}
\lim_{k\rightarrow \infty }\frac{\theta _{N+1}}{\theta _{N}}=1, %
\lim_{k\rightarrow \infty }\frac{\tau _{N+1}}{\tau _{N}}=1,  \label{tau}
\end{equation}%
where $\tau _{N}=\dim \oplus _{k=0}^{N}\mathrm{H}_{k}$. In the case of
two-point homogeneous spaces the eigenvalues $\theta _{k}$ can be written as
\[
\theta _{k}=k\left( k+\alpha +\beta +1\right) \asymp k^{2}, %
k\rightarrow \infty
\]%
if $\mathbb{M}^{d}=\mathbb{S}^{d}$, $\mathrm{P}^{d}\left( \mathbb{C}\right) $%
, $\mathrm{P}^{d}\left( \mathbb{H}\right) $, $\mathrm{P}^{16}\left( \mathrm{%
Cay}\right) $ and
\[
\theta _{2k}=2k\left( 2k+d-1\right) \asymp k^{2}, k\rightarrow \infty
\]%
if $\mathbb{M}^{d}=\mathrm{P}^{d}\left( \mathbb{R}\right) $. Note that, on
the real projective spaces, $\mathrm{P}^{d}\left( \mathbb{R}\right) $, the
only polynomials of even degree, appear because, due to the identification
of antipodal points on $\mathbb{S}^{d}$, only the even order polynomials can
be lifted to be functions on $\mathrm{P}^{d}\left( \mathbb{R}\right) $.
Applying methods of Lie algebras it is possible to show that
\[
\mathbb{S}^{d}: \alpha =\frac{d-2}{2}, \beta =\frac{d-2}{2}%
, d=2,3,\cdots ,
\]%
\[
\mathrm{P}^{d}\left( \mathbb{R}\right) : \alpha =\frac{d-2}{2},
\beta =\frac{d-2}{2}, d=2,3,\cdots ,
\]%
\[
\mathrm{P}^{d}\left( \mathbb{C}\right): \alpha =\frac{d-2}{2},
\beta =0, d=4,6,8,\cdots ,
\]%
\[
\mathrm{P}^{d}\left( \mathbb{H}\right) : \alpha =\frac{d-2}{2},
\beta =1, d=8,12,16,\cdots ,
\]%
\[
\mathrm{P}^{16}\left( \mathrm{Cay}\right) : \alpha =\frac{d-2}{2}%
, \beta =3.
\]%
It is well-known that for any two-point homogeneous space%
\[
\dim \mathrm{H}_{k}\asymp k^{d-1}
\]
and%
\[
n=\tau _{N}=\dim \oplus _{k=0}^{N}\mathrm{H}_{k}\asymp N^{d}.
\]%
See \cite{helgason}, \cite{helgason-0}, \cite{gangolli}, \cite{cartan}, \cite%
{koornwinder},\cite{bklt}, \cite{ku-to-2012}, \cite{KLT}, \cite{ku-umj},
\cite{KU-LEB-TURK}, \cite{KU-LEB-BAINOV} for more information.

To apply Theorem \ref{section-22} and Theorem \ref{theorem-33} we need
estimates of expectations $\mathbf{E}\left[ \left\Vert \cdot \right\Vert
_{B_{\left( \mathrm{J},p\right) }^{n}}\right] $. Fix $N\in \mathbb{N}$ and a
system of harmonics
\[
\left\{ \phi _{k}, 1\leq k\leq n\right\} \subset \mathcal{T}%
_{N}=\oplus _{k=0}^{N}\mathrm{H}_{k}, 
\]%
\[
\mathrm{span}\left\{ \phi _{k}, 1\leq k\leq n\right\} =\mathcal{T}_{N}%
, n=\dim \mathcal{T}_{N}.
\]%
Let \textrm{J} be the coordinate isomorphism,%
\[
\begin{array}{ccc}
\mathrm{J}:\mathbb{R}^{n} & \longrightarrow & \mathcal{T}_{N} \\
\alpha & \longmapsto & t^{\alpha }=\sum_{k=1}^{n}\alpha _{k}\phi _{k}.%
\end{array}%
\]

{\bf Lemma 3}
\label{expect} 
\begin{em}
(Lemma 3 \cite{ku-tas}) For any $2\leq p<\infty $ we have%
\begin{equation}
\mathbf{E}\left[ \left\Vert \cdot \right\Vert _{B_{\left( \mathrm{J}%
,p\right) }^{n}}\right] \leq 2^{\frac{1}{2}}\pi ^{-\frac{1}{2p}}\left(
\Gamma \left( \frac{p+1}{2}\right) \right) ^{\frac{1}{p}}.  \label{pppp11}
\end{equation}
\end{em}

Different analogues of Lemma \ref{expect} were presented in \cite{11}, \cite%
{ku-lev2}, \cite{klw}, \cite{FKK}, \cite{ku-tas}, \cite{ku-to-2014}, \cite%
{ellipsoid}.

{\bf Definition}
\begin{em}
We say that $\Lambda \in \mathcal{A}$ if $\lambda :\left( 0,\infty \right)
\mapsto \left( 0,\infty \right) $ is a decreasing continuous function such
that $\lambda \left( Ct\right) \gg \lambda \left( t\right) $, $t\rightarrow
\infty $ for any $C\geq 1$.
\end{em}

{Theorem 3}
\label{theorem-444} 
\begin{em}
Let $\mathbb{M}^{d}$ be a two-point homogeneous space
and $\Lambda =\left\{ \lambda _{k}, k\in \mathbb{N}\right\} \in
\mathcal{A}$. Then%
\[
d^{s}\left( \Lambda U_{p},L_{q}\right) \gg p^{-\frac{3}{4}}\lambda \left( s^{%
\frac{2}{d}}\right) , 1\leq q\leq p<\infty
\]%
and%
\[
d_{s}\left( \Lambda U_{p},L_{q}\right) \gg \left( \frac{q}{q-1}\right) ^{-%
\frac{3}{4}}\lambda \left( s^{\frac{2}{d}}\right) , 1<q\leq p\leq
\infty ,
\]%
as $s\rightarrow \infty $.
\end{em}

{\bf Proof}
Let $N\in \mathbb{N}$ and $n=\dim \mathcal{T}_{N}$. Since $\mathbb{M}^{d}$
is a two-point homogeneous space then the system of spherical harmonics $%
\phi _{k}$, $k\geq 0$ is essentially bounded \cite{Marzo}. Hence, applying
Theorem \ref{section-22} we get%
\begin{equation}
\mathrm{rad}\left( \mathbf{A}B_{\left( \mathrm{J},p\right) }^{n}\cap
L_{r}\left\vert B_{\left( \mathrm{J},1\right) }^{n}\right. \right) \geq
C\varrho _{n}\mathbf{E}\left[ \left\Vert \cdot \right\Vert _{B_{\left(
\mathrm{J},p\right) }^{n}}\right] ^{-\frac{3}{2}}, p\geq 2, 
\label{rad-th}
\end{equation}%
for any subspace $L_{r}\subset \mathbb{R}^{n}$, $\dim L_{r}=r$, $r\geq \frac{%
2}{3}n$. It follows from (\ref{pppp11}) that
\begin{equation}
\mathbf{E}\left[ \left\Vert \cdot \right\Vert _{B_{\left( \mathrm{J}%
,p\right) }^{n}}\right] \leq 2^{\frac{1}{2}}\left( \Gamma \left( \frac{p+1}{2%
}\right) \right) ^{\frac{1}{p}}\asymp p^{\frac{1}{2}}. \label{expp}
\end{equation}%
Also, let $\left\{ \lambda _{k}, 1\leq k\leq n\right\} $ be a
decreasing sequence and $\mathbf{A}=\mathrm{diag}\left( \lambda _{1},\cdots
,\lambda _{n}\right) $, $\det \mathbf{A}\neq 0$. Then $\left( \mathbf{A}%
^{-1}\right) ^{\ast }\mathbf{A}^{-1}$ has eigenvalues $\sigma _{k}=\lambda
_{k}^{-2}$, $1\leq k\leq n$ and therefore%
\[
\rho \left( \left( \mathbf{A}^{-1}\right) ^{\ast }\mathbf{A}^{-1}\right)
=\max \left\{ \sigma _{k}, 1\leq k\leq n\right\} =\lambda _{n}^{-2}%
,
\]%
or
\begin{equation}
\varrho _{n}=\left( \rho \left( \left( \mathbf{A}^{-1}\right) ^{\ast }%
\mathbf{A}^{-1}\right) \right) ^{-\frac{1}{2}}=\lambda _{n}.
\label{spectral}
\end{equation}%
Comparing (\ref{rad-th}), (\ref{expp}) and (\ref{spectral}) we get
\begin{equation}
\mathrm{rad}\left( \mathbf{A}B_{\left( \mathrm{J},p\right) }^{n}\cap
L_{r}\left\vert B_{\left( \mathrm{J},1\right) }^{n}\right. \right) \geq Cp^{-%
\frac{3}{4}}\lambda _{n}  \label{rad-rad}
\end{equation}%
for any subspace $L_{r}\subset \mathbb{R}^{n}$, $\dim L_{r}=r\geq \frac{2}{3}%
n$. Hence by (\ref{rad-rad}) and the definition of Gelfand widths%
\[
d^{Cn}\left( \mathbf{A}B_{\left( \mathrm{J},p\right) }^{n},\mathrm{J}\left(
L_{1}\cap \mathcal{T}_{N}\right) \right) \geq Cp^{-\frac{3}{4}}\lambda _{n}.
\]%
Consequently,
\begin{equation}
d^{n}\left( \Lambda U_{p}\cap \mathcal{T}_{N},L_{1}\right) \geq Cp^{-\frac{3%
}{4}}\lambda \left( n^{\frac{2}{d}}\right)  \label{gelf-1}
\end{equation}%
since $\Lambda \in \mathcal{A}$. Applying (\ref{dual}) we get%
\begin{equation}
d_{n}\left( \Lambda U_{\infty }\cap \mathcal{T}_{N},L_{q}\right) \geq
C\left( \frac{q}{q-1}\right) ^{-\frac{3}{4}}\lambda \left( n^{\frac{2}{d}%
}\right) ,  \label{kolm-1}
\end{equation}%
where $q>1$. Lower bounds for $d^{n}\left( \Lambda U_{p},L_{q}\right) $, $%
1\leq q\leq p<\infty $ and $d_{n}\left( \Lambda U_{p},L_{q}\right) $, $%
1<q\leq p\leq \infty $ follow from (\ref{gelf-1}) and (\ref{kolm-1}) by
imbedding. Finally, since the sequence of widths in not increasing and $%
\Lambda \in \mathcal{A}$ we get the proof for any $s\in \mathbb{N}$.
$\Box$

Observe that $\lambda \left( t\right) =t^{-\frac{\gamma }{2}}$, $\gamma >0$
in the case of Sobolev classes $W_{p}^{\gamma }$. Hence applying Theorem \ref%
{theorem-444}, (\ref{FOURIER}) and \cite{bklt} (17) p.320,
\[
\sup_{f\in W_{p}^{\gamma }}\inf_{t_{N}\in \mathcal{T}_{N}}\left\Vert
f-t_{N}\right\Vert _{q}\ll N^{-\gamma }\asymp n^{-\frac{\gamma }{d}},%
1\leq q\leq p\leq \infty ,
\]%
we get by imbedding the following statement.

{Corollary 1}
\label{cor-width} 
\begin{em}
Let $\mathbb{M}^{d}$ be a two-point homogeneous manifold.
Then for any $\gamma >0$ and $q>1$ we have
\[
d_{n}\left( W_{p}^{\gamma },L_{q}\right) \asymp n^{-\frac{\gamma }{d}}
, 1<q\leq p\leq \infty
\]%
as $n\rightarrow \infty $.
\end{em}

{Corollary 2}
\label{rim} 
\begin{em}
Let $\mathbb{M}^{d}$ be a compact, homogeneous Riemannian
manifold. Then for any $\gamma >0$ and $1<q\leq 2\leq p<\infty $ we have%
\[
\lambda ^{s}\left( W_{p}^{\gamma },L_{q}\right) \asymp s^{-\frac{\gamma }{d}},
\]%
as $s\rightarrow \infty $.
\end{em}

{\bf Proof}
Remid that $\tau _{N}=\dim \mathcal{T}_{N}$. Put $n=\tau _{N}$ and $%
\widetilde{n}=\tau _{N+1}$. Applying Weyl formula \cite{weyl}
\[
\lim_{a\rightarrow \infty }a^{-\frac{d}{2}}n\left( a\right) =\left( 2\pi ^{%
\frac{1}{2}}\right) ^{-d}\Gamma \left( 1+\frac{d}{2}\right) \mathrm{V}\left(
\mathbb{M}^{d}\right) ,
\]%
where $\mathrm{V}\left( \mathbb{M}^{d}\right) $ is the volume of $\mathbb{M}%
^{d}$ and $n\left( a\right) $ is the number of eigenvalues (each counted
with its multiplicity) smaller than $a$, we get
\[
\lim_{a\rightarrow \infty }\theta _{N}^{-\frac{d}{2}}\tau _{N}=\left( 2\pi ^{%
\frac{1}{2}}\right) ^{-d}\Gamma \left( 1+\frac{d}{2}\right) \mathrm{V}\left(
\mathbb{M}^{d}\right) .
\]%
Consequently,%
\begin{equation}
\tau _{N}\asymp \theta _{N}^{\frac{d}{2}}.  \label{tau-1}
\end{equation}%
Let $\mathbf{A}=\mathrm{diag}\left( \lambda _{1},\cdots ,\lambda _{%
\widetilde{n}}\right) $. By Theorem \ref{theorem-33}, Lemma \ref{expect} and
(\ref{tau}) we get%
\[
\mathrm{rad}\left( \mathbf{A}B_{\left( \mathrm{J},p\right) }^{\widetilde{n}%
}\cap L_{s}\left\vert B_{\left( \mathrm{J},q\right) }^{\widetilde{n}}\cap
L_{s}\right. \right) \gg \lambda \left( \theta _{N+1}\right)
\]%
\[
=\lambda \left( \frac{\theta _{N+1}}{\theta _{N}}\theta _{N}\right) \gg
\lambda \left( C\theta _{N}\right) \gg \lambda \left( \theta _{N}\right)
\]%
since $\lambda \left( t\right) =t^{-\frac{\gamma }{2}}\in \mathcal{A}$. This
means that
\[
\lambda ^{Cn}\left( W_{p}^{\gamma },L_{q}\right) \geq \lambda ^{Cn}\left(
W_{p}^{\gamma }\cap \mathcal{T}_{N+1},L_{q}\right) \gg \lambda \left( \theta
_{N}\right)
\]%
\[
\gg \lambda \left( \frac{\theta _{N}}{n^{\frac{2}{d}}}n^{\frac{2}{d}}\right)
\gg n^{-\frac{\gamma }{d}}.
\]%
where in the last line we used (\ref{tau-1}). Finally, since the seqence of
widths is nonincreasing we get the proof.
$\Box$

The author wishes to thank referees and Communicating Editor for valuable
comments which have helped clarify the presentation of the paper.

\bigskip

\end{document}